\input amstex
\documentstyle{amsppt}
\input bull-ppt
\keyedby{bull496e/paz}

\define\Ints{\Bbb Z}

\define\barr#1{\overline#1}
\predefine\prevec{\Vec}
\redefine\Vec{\hbox{\rm Vec}}

\ratitle
\topmatter
\cvol{31}
\cvolyear{1994}
\cmonth{July}
\cyear{1994}
\cvolno{1}
\cpgs{50-53}

\title 
Diffeomorphisms of Manifolds\\
with Finite 
Fundamental Group\endtitle 
\author 
 Georgia Triantafillou \endauthor 
\address Department of Mathematics, Temple University, 
Philadelphia, 
Pennsylvania 19122\endaddress 
\ml georgia\@euclid.math.temple.edu\endml
\abstract 
We show that the group ${\Cal D}(M)$ of pseudoisotopy 
classes of
diffeomorphisms of a manifold of dimension $\geq 5$ and of 
finite fundamental
group is
commensurable to an arithmetic group. As a result 
$\pi_0(\text{{\it Diff\,M}})$ 
is a group 
of finite type.\endabstract
\date September 30, 1993\enddate
\endtopmatter

\document

Let $M$ be an $n$-dimensional closed smooth manifold, 
where $n\ge 5$, and let 
{\it Diff\,M} be the group 
of diffeomorphisms of $M$. The space {\it Diff\,M} (it is 
a topological space
with the $C^\infty$-topology) is in general very
complicated and has been studied by many authors. Its 
homotopy type is
known only in some  special cases. We recall the 
well-known results
{\it Diff}$(S^2)\simeq O(3)$ by Smale [Sm] and 
{\it Diff}$(S^3)\simeq O(4)$ by Hatcher [H].

The component group $\pi_0(\text{{\it Diff\,M}})$ has also 
been initially
computed in 
special cases, for instance for spheres, where the group
$\pi_0(\text{{\it Diff}}_+\,S^n)$ of 
orientation-preserving diffeomorphisms
is isomorphic to the group of homotopy spheres of 
dimension $n+1$ for $n>4$
[KM], and for products of spheres [B, Tu].
The first general result about the nature of 
$\pi_0(\text{{\it Diff\,M}})$
 for simply 
connected manifolds $M$ is due to
Sullivan [S], who showed that if $M$ is a smooth closed 
orientable {\it simply 
connected}
manifold of dimension $\geq 5$, then $\pi_0(\text{{\it 
Diff\,M}})$
 is commensurable 
to an arithmetic group.

Arithmetic groups and their properties have been studied 
in [BH] by Borel 
and Harish-Chandra. In particular they showed that every 
arithmetic 
group is finitely presented. It follows that 
$\pi_0(\text{{\it Diff\,M}})$
is finitely presented if $M$ is simply connected.
In fact $\pi_0(\text{{\it Diff\,M}})$ in this case is a 
group of finite type by a 
result of Borel and Serre that arithmetic groups are of 
finite
type [BS]. By definition a group $\pi$ is of {\it finite 
type} if  its 
classifying space $B\pi$ is homotopy equivalent to a 
CW-complex with
finitely many cells in each dimension. Being of finite 
type for a group
implies and in fact is much stronger than finite 
presentation. 
Two groups are said
to be {\it commensurable} if there is a finite
sequence of homomorphisms between them
which have finite kernels and images of finite index. 

On the other hand $\pi_0(\text{{\it Diff\,M}})$
 can be very large if $M$ is not simply 
connected. For instance,
 $\pi_0\text{({\it Diff}}\,T^n)$ contains a 
subgroup isomorphic to
a direct sum of infinitely many copies of ${\Ints}_2$ for 
$n\ge 5$, and 
$\pi_0(\hbox{\it Diff}(T^n\times S^2))$
contains a free abelian subgroup of infinite 
rank for $n\ge 3$, where $T^n$ is the torus [HS]. In these 
cases $\pi_1(M)$ is
infinite. 


In this paper we study $\pi_0(\text{{\it Diff\,M}})$ and 
${\Cal D}(M)$
in the case where $\pi_1(M)$
is finite, where ${\Cal D}(M)$ is the group of 
pseudoisotopy classes of 
diffeomorphisms of $M$. 

Pseudoisotopy is defined by dropping the level
preserving requirement of isotopy, namely, two 
diffeomorphisms $f$ and $g$ of
$M$ are pseudoisotopic if there is a diffeomorphism 
$F:M\times I\to M\times
I$ which restricts to $f$ and $g$ at the two ends of 
$M\times I$ respectively.
While in the simply connected case $\pi_0(\hbox{\it 
Diff\,M})={\Cal D}(M)$ if 
$n\geq 5$ by Cerf's result [C], in the nonsimply connected 
case there is 
in general an Abelian kernel of the forgetful map
$$0\to A\to \pi_0(\hbox{\it Diff\,M})\to {\Cal D}(M)\to 0,$$
where $A$ has been computed by Hatcher and Wagoner [HW] 
and Igusa [I]. It is also
well known to the experts that $A$ is finitely generated 
if $\pi_1(M)$ is 
finite. We show:

\proclaim{ Theorem 1} Let $M$ be a closed, smooth, 
orientable manifold
of dimension $\ge 5$ with 
finite
fundamental group. Then the group ${\Cal D}(M)$ of 
pseudoisotopy classes of 
diffeomorphisms of $M$ is commensurable to an arithmetic 
group.
\endproclaim

As a consequence of this and the exact sequence above we 
get:

\proclaim{ Theorem 2 }If $M$ is a smooth, closed, 
orientable manifold of dimension
$\geq 5$ with finite fundamental group, then 
$\pi_0(\text{{\it Diff\,M}})$ is of finite 
type.
\endproclaim

\heading Outline of proofs \endheading

The proof breaks up naturally into two parts: a homotopy 
theoretical part 
and a geometric part. Let $aut(M)$ be the group of 
homotopy classes of 
self-homotopy equivalences of $M$, and let 
$aut_t(M)\subseteq aut(M)$ be the
subgroup, the elements of which are represented by maps 
$f$ which preserve
the tangent bundle, i.e., $f^*(T(M))=T(M)$.
First we show that the group $aut_t(M)$ is commensurable 
to an arithmetic 
group. For this purpose we employ minimal model techniques 
in the equivariant 
context and certain facts from the theory of algebraic 
and arithmetic groups. In the presence of a 
finite fundamental group we work with the minimal model of 
the universal
cover of $M$, where the model is now equipped with an 
action of $\pi_1$.
By using this equivariant minimal model, we show the 
following. 


\proclaim{ Theorem 3} The group $aut(M)$ of homotopy 
classes of
self-homotopy equivalences of a finite CW-complex $M$ with 
finite 
fundamental group is commensurable to an arithmetic group. 
\endproclaim


This generalizes a result of 
[W] for the simply connected case and independently of
[S] for the nilpotent case. It has been shown by 
[DDK]  that $aut(X)$ is of finite type for virtually 
nilpotent spaces 
including spaces with finite fundamental group. The full 
strength
of arithmeticity of $aut(X)$ above is needed to prove the 
following
result.
We note here that the properties of being of finite type, 
finitely 
presented, finitely generated, or arithmetic are not in 
general inherited
by subgroups.

\proclaim{ Theorem 4} The group $aut_t(M)$ of homotopy 
classes of
tangential self-homotopy equivalences of a manifold $M$ 
with finite 
fundamental group is commensurable to an arithmetic group. 
\endproclaim


In the geometric part of the proof we 
apply  surgery  techniques and the Atiyah-Singer
$G$-signature theorem to relate 
the group of tangential homotopy equivalences to  
the group {\it diff\,M} of homotopy classes of 
diffeomorphisms, where the 
latter is the
image of $\pi_0(\text{{\it Diff\,M}})$ into the group of 
homotopy classes of 
self-homotopy equivalences. Out of this comparison follows 
in particular:
	
\proclaim{ Theorem 5} The group of homotopy classes of 
tangential 
simple self-homotopy equivalences $aut_t^s(M)$
is commensurable to an arithmetic group if $M$ has finite 
fundamental group
and dimension $\ge5$.
\endproclaim


Further, we relate bundle homotopy classes of tangential
homotopy equivalences to pseudoisotopy classes of 
diffeomorphisms.
Let $Aut(M, t_0)$ be the set (monoid) of self-homotopy 
equivalences of $M$
which preserve the tangent bundle except at a  base
point $x_0\in M$. More precisely, the 
elements of $Aut(M, t_0)$ are  
pairs $(h, b)$, where $h$ is a   base point preserving
self-homotopy equivalence of $M$ and $b:T(M-{x_0})\to 
T(M-{x_0})$ is a 
bundle map covering $h$. Let $aut(M, t_0)=\pi_0(Aut(M, 
t_0))$
be the group of bundle homotopy classes of such pairs.

The following is the main geometric result leading to the 
arithmeticity of
${\Cal D}(M)$ as it relates the latter to homotopy and 
hence algebraic data.
Under the same conditions as above we have:

\proclaim{ Theorem 6} The group of pseudoisotopy classes 
of diffeomorphisms 
${\Cal D}(M)$ is commensurable to the group $aut(M, t_0)$ 
of bundle homotopy
classes of self-homotopy equivalences  covered by bundle 
maps of the
tangent bundle except at one point.
\endproclaim

These groups appear as  arithmetic subgroups of an algebraic
matrix group over the rationals involving isomorphisms of 
the minimal 
model provided with certain additional structure. 

Details of the proofs will appear elsewhere.
\heading Acknowledgments\endheading

I would like to thank Shmuel Weinberger for useful 
discussions and the 
University of Chicago for its hospitality while part of this
research was done.

\enddocument